\documentclass[a4paper,10pt]{article}
\usepackage{amsmath,amssymb}
\usepackage[latin1]{inputenc}
\usepackage[dvips]{epsfig}
\usepackage{mathrsfs}

\newcommand{\TA}{\mathbb{A}}
\newcommand{\TD}{\mathbb{D}}
\newcommand{\TE}{\mathbb{E}}

\newcommand{\pl}{\operatorname{PreLie}}
\newcommand{\dend}{\operatorname{Dend}}
\newcommand{\mld}{\operatorname{Mould}}
\newcommand{\ass}{\operatorname{Assoc}}
\newcommand{\NCP}{\operatorname{NCP}}
\newcommand{\tri}{\operatorname{TriDend}}
\newcommand{\udim}{\underline{\dim}}

\newcommand{\sh}{\operatorname{Sh}}

\newcommand{\sym}{\mathfrak{S}}
\newcommand{\one}{\mathsf{\bf{1}}}
\newcommand{\res}{\operatorname{Res}}

\newcommand{\bY}{\mathcal{Y}}
\newcommand{\forg}{\mathscr{F}}
\newcommand{\sous}{\curvearrowleft}

\newcommand{\lft}{\epsfig{file=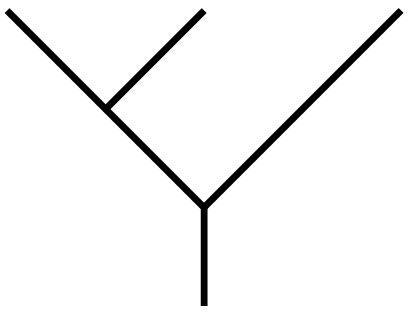,height=2mm}}
\newcommand{\rgt}{\epsfig{file=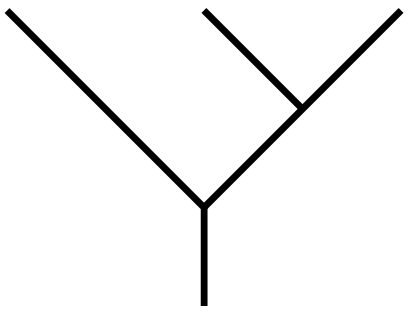,height=2mm}}
\newcommand{\mil}{\epsfig{file=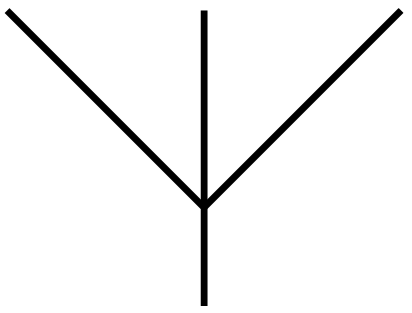,height=2mm}}
\newcommand{\gA}{\epsfig{file=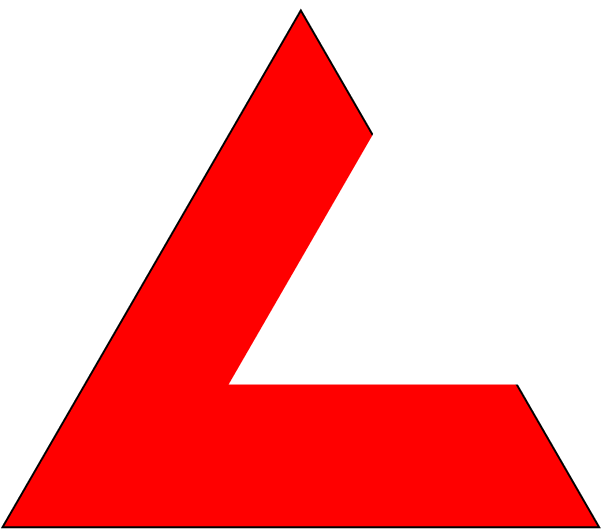,height=3mm}}
\newcommand{\dA}{\epsfig{file=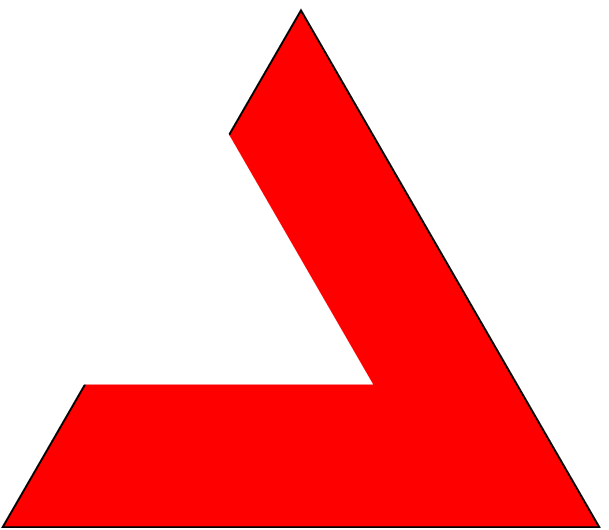,height=3mm}}
\newcommand{\mA}{\epsfig{file=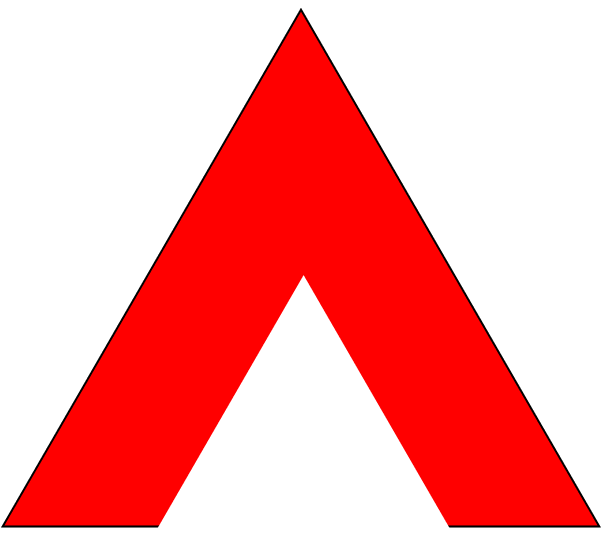,height=3mm}}

\newtheorem{theorem}{Theorem}[section] 
\newtheorem{proposition}[theorem]{Proposition} 
\newtheorem{conjecture}[theorem]{Conjecture} 
\newtheorem{corollary}[theorem]{Corollary} 
\newtheorem{lemma}[theorem]{Lemma}

\newtheorem{remark}[theorem]{Remark}

\newenvironment{proof}{\begin{trivlist}\item{\bf{Proof.}}}
  {\hfill\rule{2mm}{2mm}\end{trivlist}}

\title{The anticyclic operad of moulds} 
\author{F. Chapoton}
\date{\today}

\begin{document}

\maketitle

\begin{abstract}
  A new anticyclic operad $\mld$ is introduced, on spaces of functions
  in several variables. It is proved that the Dendriform operad is an
  anticyclic suboperad of this operad. Many operations on the free
  $\mld$ algebra on one generator are introduced and studied. Under
  some restrictions, a forgetful map from moulds to formal vector
  fields is then defined. A connection to the theory of tilting
  modules for quivers of type $\TA$ is also described.
\end{abstract}

\setcounter{section}{-1}

\section{Introduction}
The aim of this article is to build and use a new connection between
the theory of operads and the theory of moulds. Operads were
introduced in algebraic topology in the 1960's. After being somewhat
neglected for some decades, this notion has found a new impetus
recently, in connection with mathematical physics, moduli spaces of
curves and algebraic combinatorics.

Moulds have a rather different origin. They have been introduced in
analysis by J. Ecalle, as a convenient tool to handle complicated
singular functions, in relation with his theory of resurgence. Later,
he developed around moulds a large apparatus which allowed him to make
substantial progress in the theory of polyzetas
\cite{eca2002,eca2003,eca2004}. In this article, only the simplest
case of moulds will be considered, not the more general case of
bimoulds.

The first result of this article is the existence of a very simple
structure of operad on moulds, denoted by $\mld$. In fact, one can
define on moulds the finer structure of an anticyclic operad,
involving in addition to composition maps some actions of cyclic
groups.

This structure is then shown to contain as an anticyclic suboperad the
so-called Dendriform operad introduced by J.-L. Loday \cite{livre} and
denoted by $\dend$, which has been much studied recently
\cite{planarBT,lodayronco}. This provides a radically new point of
view on the operad $\dend$ and the key to some new results. The main
result is the explicit description of the smallest subset of $\dend$
containing its usual generators and closed under the anticyclic operad
structure, by the mean of new combinatorial objects called
non-crossing plants.

This article also contains the description of many different
operations on moulds, coming either from the operad or from the mould
viewpoint, and some of their properties. This is used to prove the
existence of a morphism from the Lie algebra of moulds (for one of the
Lie brackets and under some restrictions) to the Lie algebra of formal
vector fields in one indeterminate. Interesting and natural examples
of moulds are provided and their images by this map are computed.

In some sense, this article provides a reformulation of the basic
results of Ecalle on moulds in a more classical algebraic language.
This includes notably the so-called $\mathrm{ARI}$ bracket, for which
we provide a very short definition using the operations obtained from
the operad structure. We use this definition to prove some properties
of this bracket. It should be said that our setting does not seem to
extend to bimoulds, hence can only describe a small part of the theory
of Ecalle.

In the last section, it is recalled that the Dendriform operad is
strongly related to the theory of tilting modules for the
equi-oriented quivers of type $\TA$, and how the results of the
present article fit very-well in this relationship. Some conjectural
extension of the properties in type $\TA$ to other Dynkin diagrams are
proposed.

Many thanks to Jean Ecalle for discussions which have led to this
research.

MSC 2000: 18D50, 05C05, 05E

\section{Notations and definitions}

\label{rappels}

We recall here the terminology we need concerning moulds.

A \textit{mould} is a sequence $(f_n)_{n\geq 1}$, where $f_n$ is a
function of the variables $\{u_1,\dots,u_n\}$. A mould is said to have
degree $n$ if its only non-zero component is $f_n$. In this case, this
unique component will be denoted $f$ by a convenient abuse of notation.

A mould $f$ of degree $n$ is called \textit{alternal} if it satisfies the
following conditions, for $1\leq i \leq n-1$:
\begin{equation}
  \sum_{\sigma\in\sh(i,n-i)} f(u_{\sigma(1)},\dots,u_{\sigma(n)})=0,
\end{equation}
where $\sigma$ runs over the shuffle permutations of
$\{u_1,\dots,u_i\}$ and $\{u_{i+1},\dots,u_n\}$, i.e. permutations
such that $\sigma(1)<\dots<\sigma(i)$ and
$\sigma(i+1)<\dots<\sigma(n)$. One can then extend this definition: a
mould is called alternal if each of its components is alternal.

A mould $f$ of degree $n$ is called \textit{vegetal} if it satisfies the
following equation:
\begin{equation}
  {u_1\dots u_n} \sum_{\sigma\in\sym_n} f(t\,
  u_{\sigma(1)},\dots,t\, u_{\sigma(n)})=n! f(t,\dots,t),
\end{equation}
where $\sym_n$ is the permutation group of $\{1,\dots,n\}$.

There is a natural associative product on moulds, defined for $f$ of
degree $m$ and $g$ of degree $n$ by
\begin{equation}
  \label{def_MU}
  \mathrm{MU}(f,g)=f(u_1,\dots,u_m)g(u_{m+1},\dots,u_{m+n}).
\end{equation}

The associated Lie bracket is
\begin{equation}
  \label{def_LIMU}
  \mathrm{LIMU}(f,g)=\mathrm{MU}(f,g)-\mathrm{MU}(g,f).
\end{equation}

We will often use the following convenient shorthand notation:
\begin{equation}
  u_{i..j}:=\sum_{i \leq k \leq j} u_k.
\end{equation}

At some places in the text, we will use the following shorthand
notation. For a shuffle $\sigma$ of two ordered sets $S'$ and $S''$,
let $u_\sigma$ be the sequence $u_s$ for $s \in S' \cup S''$ in the
order specified by $\sigma$.

\section{The Mould operad}

For $n \geq 1$, let $\mld(n)$ be the vector space of rational
functions with rational coefficients in the variables
$\{u_1,\dots,u_n\}$. We will show that the collection
$\mld=(\mld(n))_{n\geq 1}$ has the structure of an anticyclic non-symmetric
operad. The reader is referred to \cite{markl,markl_livre} for the
basics of the theory of operads and anticyclic operads.

First, let $\one$ be the function $1/u_1$ in $\mld(1)$. This will be
the unit of the operad.

Let us then introduce a map $\tau$, which is called the \textit{push}.
It is defined on $\mld(n)$ by
\begin{equation}
   \tau(f)(u_1,\dots,u_n)=f(-u_{1\dots n},u_1,\dots,u_{n-1}).
\end{equation}
Note that $\tau$ has order $n+1$ on $\mld(n)$. It will give the cyclic
action of the operad. Let us note also that $\tau(\one)=-\one$. This is one of the axioms of an
anticyclic operad.

Let us now introduce the composition maps $\circ_i$ from
$\mld(m)\otimes \mld(n)$ to $\mld(m+n-1)$, with $1 \leq i \leq m$. Let
$f$ be in $\mld(m)$ and $g$ be in $\mld(n)$. The function $f \circ_i
g$ is defined by
\begin{equation}
  \label{main}
  u_{i\dots{i+n-1}} f(u_1,\dots,u_{i-1},u_{i\dots{i+n-1}},u_{i+n},\dots,u_{m+n-1}) g(u_i,\dots,u_{i+n-1}).
\end{equation}

\begin{theorem}
  The push $\tau$ and composition maps $\circ_i$ define the structure
  of an anticyclic non-symmetric operad $\mld$.
\end{theorem}

\begin{proof}
  One has first to check that these composition maps do indeed define
  a non-symmetric operad. The unit $\one$ has clearly the expected
  properties: $\one \circ_1 f=f$ and $f\circ_i \one=f$ for all
  $f\in\mld(m)$ and $1\leq i\leq m$.  One has also to check two
  ``associativity'' axioms.

  Let $f,g,h$ be in $\mld(m),\mld(n)$ and $\mld(p)$.

  Let $i$ and $j$ be such that $1 \leq i <j \leq m$. Then one has to check that
  \begin{equation}
    (f \circ_i g) \circ_{j+n-1} h = (f\circ_j h)\circ_i g.
  \end{equation}

  Indeed, both sides are equal to
  \begin{multline}
    u_{i\dots{i+n-1}}\, u_{{j+n-1}\dots{j+n+p-2}}\,
    g(u_i,\dots,u_{i+n-1})h(u_{j+n-1},\dots,u_{j+n+p-2})\\f(u_1,\dots,u_{i-1},u_{i\dots{i+n-1}},u_{i+n},\\
    \dots,u_{j+n-2},u_{{j+n-1}\dots{j+n+p-2}},u_{j+p+n-1},\dots,u_{m+n+p-2}).
  \end{multline}

  Let now $i$ and $j$ be such that $1 \leq i \leq m$ and $1\leq j \leq
  n$. One has to check that
  \begin{equation}
    f \circ_i (g \circ_j h) = (f \circ_i g ) \circ_{j+i-1} h. 
  \end{equation}
  Indeed, both sides are equal to
  \begin{multline}
    f(u_1,\dots,u_{i-1},u_{i\dots{i+p+n-2}},u_{i+p+n-1},\dots,u_{m+n+p-2})\\
    g(u_i,\dots,u_{j+i-2},u_{{j+i-1}\dots{j+i+p-2}},u_{j+i+p-1},\dots,u_{i+p+n-2})
    \\ u_{i\dots{i+p+n-2}}\, u_{{j+i-1}\dots{j+i+p-2}} h(u_{j+i-1},\dots,u_{j+i+p-2}).
  \end{multline}
  
  This proves that $\mld$ is a non-symmetric operad. Then one has to
  verify that $\tau$ gives furthermore an anticyclic structure on this
  operad.
  
  For this, one has to check two identities. The first one is
  \begin{equation}
    \tau(f\circ_i g)=\tau(f) \circ_{i-1} g,
  \end{equation}
  for $f\in \mld(m)$, $g\in\mld(n)$ and $2 \leq i\leq m$. Indeed, this holds true, as both sides are equal to
  \begin{multline}
    f(-u_{1\dots m+n-1},u_1,\dots,u_{i-2},u_{i-1\dots
      n+i-2},u_{n+i-1},\dots,u_{m+n-2}) \\ g(u_{i-1},\dots,u_{n+i-2})
    u_{i-1\dots n+i-2}.
  \end{multline}
  
  The other identity that we have to check is
  \begin{equation}
    \tau(f \circ_1 g)=-\tau(g) \circ_n \tau(f),
  \end{equation}
  for $f\in \mld(m)$ and $g\in\mld(n)$. Again, this is true as both
  sides are equal to
  \begin{equation}
    -u_{n\dots m+n-1} f(-u_{n\dots m+n-1},u_n,\dots,u_{m+n-2}) g(-u_{1\dots m+n-1},u_1,\dots,u_{n-1}).
  \end{equation}
  
\end{proof}

\begin{remark}\rm
  \label{homogen}
  It follows from Eq. (\ref{main}) that if $f$ and $g$ are homogeneous
  functions of weight $d_f$ and $d_g$ (where all variables $u_i$ are
  taken with weight $1$), then $f\circ_i g$ is also homogeneous of
  weight $d_f+d_g+1$. Also, the action of $\tau$ clearly preserves the
  weight in that sense. Hence the collection of subspaces of
  homogeneous rational functions of weight $-n$ in $\mld(n)$ is a
  anticyclic suboperad.
\end{remark}

\begin{remark}\rm
   Another consequence of Eq. (\ref{main}) is the following. Let $H_n$
   be the product
   \begin{equation}
     H_n=\prod_{1\leq i \leq j \leq n} u_{i\dots j}.
   \end{equation}
   One can see that, if $H_m f$ and $H_n g$ are polynomials, then so
   is $H_{m+n-1} (f\circ_i g)$. It is also true that $H_m \tau(f)$ is
   polynomial if $H_m f$ is, as one can easily check that $\tau$
   preserves $H_m$ up to sign.
\end{remark}

\smallskip

Combining the two previous remarks, one gets that the subspace of
$\mld(n)$ made of homogeneous rational functions $f$ of weight $-n$
such that $H_n f$ is a polynomial define a anticyclic suboperad,
which is finite dimensional in each degree.

\begin{remark}\rm
  \label{nicepole}
  By a similar argument, one can also note that the subspace of
  rational functions that have only poles of the shape $u_{i\dots j}$
  (at some power) for some $i\leq j$ is also stable for the
  composition and the cyclic action. Such functions will be said to
  have \textit{nice poles}.
\end{remark}

\section{The Dendriform operad}

\label{dendri_op}

The Dendriform operad was introduced by Loday \cite{livre}, motivated
by some problem in algebraic topology. Later, it was shown to be an
anticyclic operad \cite{anticyclic}. We refer the reader to the book
\cite{livre} for more details on this operad.

Recall that the Dendriform operad is an operad in the category of
vector spaces, generated by $\lft$ and $\rgt$ of degree $2$ with
relations
\begin{align}
  \lft \circ_1 (\lft + \rgt)&=\lft \circ_2 \lft,\\
  \label{dend}
  \rgt \circ_1 \lft &= \lft \circ_2 \rgt,\\
  \rgt \circ_1 \rgt &= \rgt \circ_2 (\lft + \rgt).
\end{align}

The dimension of $\dend(n)$ is the Catalan number
\begin{equation}
  c_n = \frac{1}{n+1}\binom{2n}{n+1}.
\end{equation}

There is a basis of $\dend(n)$ indexed by the set $\bY(n)$ of rooted
planar binary trees with $n+1$ leaves. In the presentation above,
$\lft$ and $\rgt$ correspond to the two planar binary trees in
$\bY(2)$. We will sometimes denote by $Y$ the unique planar binary
tree of degree $1$. 

\begin{figure}
  \begin{center}
    \epsfig{file=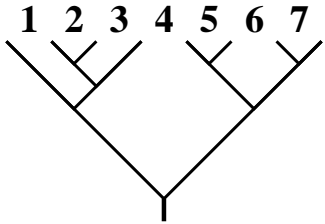,height=3cm} 
    \caption{a planar binary tree $T$ and the standard numbering}
    \label{planar}
  \end{center}
\end{figure}

The cyclic action is defined on the generators by
\begin{align}
  \tau(\lft)&=\rgt,\\
  \tau(\rgt)&=-(\lft+\rgt).
\end{align}

\begin{theorem}
  There is a unique map $\psi$ of anticyclic non-symmetric operads from
  $\dend$ to $\mld$ which maps $\lft$ to $1/(u_1 u_{1\dots 2})$ and $\rgt$
  to $1/(u_{1\dots 2} u_2)$.
\end{theorem}

\begin{proof}
  It is quite immediate to check, using the known quadratic binary
  presentation of $\dend$ and the description of the cyclic action on
  these generators recalled above, that this indeed defines a morphism
  $\psi$ of operads and that this morphism $\psi$ is a morphism of
  anticyclic operads.
\end{proof}

Let us now describe the image by $\psi$ of a planar binary tree $T$ in
$\bY(n)$.

Let us define, for each inner vertex $v$ of $T$, a linear function
$\udim(v)$ in the variables $u_1,\dots,u_n$. One can label from left
to right the spaces between the leaves from $1$ to $n$ as in Fig.
\ref{planar}. Then the vertex $v$ defines a pair of leaves (its
leftmost and rightmost descendants), enclosing a subinterval $[i,j]$
of $[1,n]$. Let $\udim(v)$ be $u_{i\dots j}$.

\begin{proposition}
  \label{fraction_arbre}
  Let $T$ be a planar binary tree. Then its image $\psi(T)$ is the
  inverse of the product of factors $\udim(v)$ over all inner vertices
  $v$ of $T$.
\end{proposition}

For instance, the image of the tree $T$ of Fig. \ref{planar} is
\begin{equation}
  \psi(T)= 
  1/\left(u_{1\dots 3} u_2 u_{2\dots 3}u_{1\dots 7} u_5 u_{5\dots 7} u_7  \right).
\end{equation}

\begin{proof}
The proof is by induction on $n$. The proposition is true for $n=1$ or
$2$.

Assume that the Proposition is true up to degree $n$. Let $T$ be a
planar binary tree in $\bY(n+1)$. By picking a top vertex of $T$ (any
inner vertex of maximal height), one can find a tree $S$ in $\bY(n)$
and an index $i$ such that $T=S \circ_i \lft$ or $T=S \circ_i \rgt$.

Then one can check that the description given above has the correct
behavior with respect to such compositions in $\mld$ and in $\dend$.
\end{proof}

Let us now introduce a classical map $\pi$ from permutations of
$\{1,\dots,n\}$ to planar binary trees in $\bY(n)$. First, note that
one can use the standard numbering as in Fig. \ref{planar} to label
the inner vertices of a planar binary tree from $1$ to $n$ from left
to right. Then each tree $T$ induces a natural partial order $\leq_T$
on $\{1,\dots,n\}$ by saying that $i \leq_T j$ if the inner vertex $i$
is below (i.e. is an ancestor of) the inner vertex $j$. The map $\pi$
is characterized by the property that $\pi(\sigma)=T$ if and only if
the total order $\sigma(1)<\sigma(2)<\dots<\sigma(n)$ is an extension
of the partial order $\leq_T$. In particular, $\sigma(1)$ must be the
index of the bottom inner vertex of $T$. The map $\pi$ is surjective
and has a standard construction by induction, see for example
\cite{planarBT}. For instance, the images by $\pi$ of the permutations
$4163527$ and $4651372$ are both the tree of Fig. \ref{planar}.

Let us define the multi-residue of an element in $\mld(n)$ according to
a permutation $\sigma \in \sym_n$:
\begin{equation}
  \oint_\sigma f = (2 i \pi)^{-n} \oint_{u_{\sigma(1)}} \dots \oint_{u_{\sigma(n)}} f.
\end{equation}

\begin{proposition}
  For a planar binary tree $T$ in $\bY(n)$ and a permutation $\sigma
  \in \sym_n$, the multi-residue $\oint_\sigma \psi(T)$ does not
  vanish if and only if $\pi(\sigma)=T$.
\end{proposition}
\begin{proof}
  The proof is by induction. The statement is clear if $n=1$. Let $k$
  be $\sigma(n)$.

  Let us assume that $\pi(\sigma)=T$. By the discussion above on the
  properties of $\pi$, this implies that the vertex $k$ is a top
  vertex of $T$ (a maximal element for $\leq_T$). Then by computing
  the innermost residue with respect to $u_k$, the multi-residue
  $\oint_\sigma \psi(T)$ reduces to the multi-residue of the function
  obtained by replacing $u_k$ by $0$ in $u_k \psi(T)$, with respect to
  indices $\{1,\dots,n\}$ but $k$, in some order. By renumbering the
  variables, this multi-residue is just $\oint_{\sigma'} \psi(T')$
  where $T'$ is obtained by removing the top vertex $k$ from $T$ and
  $\sigma'$ is the induced permutation of $\{1,\dots,n-1\}$. It is
  clear that $\pi(\sigma')=T'$, hence the residue $\oint_{\sigma'}
  \psi(T')$ is not zero by induction and therefore the residue
  $\oint_\sigma \psi(T)$ is not zero too.

  Let us now assume that $\pi(\sigma)\not=T$. If the vertex numbered
  $k$ is not a top vertex , then the
  residue with respect to $u_k$ is zero, as $u_k$ is not a pole of
  $\psi(T)$.  If the vertex number $k$ is a top vertex, then, with the
  same notations as above, one necessarily has that
  $\pi(\sigma')\not=T'$. Hence by induction $\oint_{\sigma'}
  \psi(T')=0$ and therefore $\oint_\sigma \psi(T)=0$.

\end{proof}

\begin{theorem}
  The morphism $\psi$ is injective. 
\end{theorem}
\begin{proof}
  It is enough to prove that the the functions $\psi(T)$ for all planar
  binary trees in $\bY(n)$ are linearly independent. This follows
  from the previous proposition.   
\end{proof}
  
\subsection{The Associative operad}

It is known that the Associative operad is the suboperad of the
Dendriform operad generated by $\lft+\rgt$. Furthermore, the basis of
the one-dimensional space $\ass(n)$ is mapped in $\dend$ to the sum of
all planar binary trees.

Inside the Mould operad, as the image of $\lft+\rgt$ is $1/(u_1 u_2)$,
one can check by induction that the image of the basis of $\ass(n)$ is
the inverse of the product $u_1 \dots u_n$.

\begin{remark}\rm
  the Associative operad is not stable for the anticyclic structure,
though. The smallest anticyclic suboperad of $\mld$ containing $\ass$
is $\dend$.
\end{remark}

\subsection{The graded Tridendriform operad}

Several generalizations of the Dendriform operad have been introduced
in \cite{Koperade} and \cite{loday_tri}. They all share the common
point that they have in degree $n$ a basis indexed by all planar trees
instead of just planar binary trees.

Let us consider among them the operad $gr\tri$ which is the associated
graded operad of the Tridendriform operad (which is a filtered
operad). It has been considered both in \cite{Koperade} and
\cite{loday_tri} and contains as a suboperad the Dendriform operad.
The operad $gr\tri$ is generated by the dendriform generators together
with another associative operation $\mil$ in degree $2$. One can
deduce from the results of \cite{loday_tri} a presentation by
generators and relations of $gr\tri$. It consists of the associativity
relations for $\mil$, of the $3$ relations for the dendriform
generators $\lft$ and $\rgt$ given at the beginning of \S
\ref{dendri_op}, and of $3$ more relations:
\begin{align}
  \mil \circ_1 \rgt &= \mil \circ_2 \lft,\\
  \mil \circ_1 \lft &= \lft \circ_2 \mil,\\
  \rgt \circ_1 \mil &= \mil \circ_2 \rgt.
\end{align}

\begin{proposition}
  By extending the morphism $\psi$ by $\mil \mapsto
  \frac{1}{u_1+u_2}$, one gets a morphism (still denoted by $\psi$) of
  operads from $gr\tri$ to $\mld$.
\end{proposition}

\begin{proof}
  This is an immediate verification.
\end{proof}

We will not check here whether or not this morphism is injective. This
may follow from the same kind of arguments as for the Dendriform
operad. Indeed, the image of a planar tree by $\psi$ has a simple
description given by an obvious extension of Prop.
\ref{fraction_arbre}.

\begin{remark}
  There is no known cyclic or anticyclic structure on $gr\tri$. The
  image of $gr\tri$ by $\psi$ is not closed for the action of
  $\tau$. One may ask for a description of the closure of $gr\tri$ in
  the anticyclic operad $\mld$. 
\end{remark}

\section{Free algebra on one generator}

Let us consider the free algebra on one generator for the operad
$\mld$. This can be identified with the direct sum of all spaces
$\mld(n)$, which will also be denoted by $\mld$.

\subsection{Dendriform products}

The inclusion of the operad $\dend$ in $\mld$ defines the structure
of a dendriform algebra on the free $\mld$ algebra on one generator:
we have two binary operations $\succ$ and $\prec$ defined for
$f\in\mld(m)$ and $g\in\mld(n)$ by
\begin{equation}
  f \succ g = \left(\frac{1}{u_1 u_{1\dots 2}}\circ_2 g\right) \circ_1
  f = f(u_1,\dots,u_m)g(u_{m+1},\dots,u_{m+n})\frac{u_{m+1\dots
  m+n}}{u_{1\dots m+n}}
\end{equation}
and
\begin{equation}
  f \prec g = \left(\frac{1}{u_{1\dots 2} u_2}\circ_2 g\right) \circ_1
  f= f(u_1,\dots,u_m)g(u_{m+1},\dots,u_{m+n})\frac{u_{1\dots
  m}}{u_{1\dots m+n}}.
\end{equation}

\begin{proposition}
  \label{dend_vegetal}
  If $f$ and $g$ are vegetal, then so are $f\succ g$ and $f\prec g$.
\end{proposition}

\begin{proof}
  Let us do the proof for $f \succ g$, the other case being similar.
  One has to compute
  \begin{multline}
    \frac{u_1 \dots u_{m+n}}{u_{1\dots
        m+n}}\sum_{\sigma\in\sym_{m+n}} f(t\,
    u_{\sigma(1)},\dots,t\,u_{\sigma(m)})g(t\,u_{\sigma(m+1)},\dots,t\,u_{\sigma(m+n)})\\(u_{\sigma(m+1)}+\dots+u_{\sigma(m+n)}).
  \end{multline}
  Let us introduce the set $E=\{{\sigma(m+1)},\dots,{\sigma(m+n)}\}$.
  Then one can rewrite the previous sum as
  \begin{multline}
    \frac{u_1 \dots u_{m+n}}{u_{1\dots m+n}}\sum_E \left({\sum_{i\in
        E}u_i} \right) \sum_{\sigma}
    f(t\,u_{\sigma'(1)},\dots,t\,u_{\sigma'(m)})\\\sum_{\sigma''}
    g(t\,u_{\sigma''(m+1)},\dots,t\,u_{\sigma''(m+n)}) ,
  \end{multline}
  where $E$ runs over the set of subsets of cardinal $n$ of
  $\{1,\dots,m+n\}$, $\sigma'$ is a bijection from $\{1,\dots,n\}$ to
  the complement of $E$ and $\sigma''$ is a bijection from
  $\{m+1,\dots,m+n\}$ to $E$.  Using the vegetal property of $f$ and
  $g$, this reduces to
  \begin{equation}
    \frac{m! n!}{u_{1\dots      m+n}}f(t,\dots, t)
      g(t,\dots, t) \sum_E 
      {\sum_{i\in E}u_i}.
  \end{equation}
   Reversing summations, this gives
  \begin{equation}
    {(m+n-1)! n} f(t\,\dots, t) g(t\,\dots, t),
  \end{equation}
  which is $(m+n)! (f\succ g)(t,\dots,t)$. Hence $f\succ g$ is vegetal.
\end{proof}

\subsection{Associative product}

The inclusion of the Associative operad in the Mould operad implies
that the formula
\begin{equation}
  \mathrm{MU}(f,g)=\left(\frac{1}{u_1 u_2}\circ_2 g\right) \circ_1 f
\end{equation}
defines an associative product on the free $\mld$-algebra on one
generator. One can check that this associative product is exactly the
product called $\mathrm{MU}$ in the terminology of moulds, see (\ref{def_MU}).
Hence the associated bracket is the so-called $\mathrm{LIMU}$ bracket. Note also
that $f\succ g+f\prec g=\mathrm{MU}(f,g)$.

As a consequence of Prop. \ref{dend_vegetal}, one has
\begin{corollary}
  \label{mu_vegetal}
  If $f$ and $g$ are vegetal, then so are $\mathrm{MU}(f,g)$ and
  $\mathrm{LIMU}(f,g)$.
\end{corollary}

\subsection{Pre-Lie product}

Recall that a pre-Lie product on a vector space $V$ is a bilinear map
$\sous$ from $V$ to $V$ such that 
\begin{equation}
  (x \sous y) \sous z-x \sous (y \sous z)
=(x \sous z) \sous y-x \sous (z \sous y).
\end{equation}
This notion is related to manifolds with affine structures and to
groups with left-invariant affine structures. As for associative
algebras, the corresponding antisymmetric bracket $[x,y]=x\sous y-y
\sous x$ is a Lie bracket.  For a reference on pre-Lie algebras, the
reader may consult \cite{prelie}.

As there is a injective morphism from the $\pl$ operad to
the symmetric version of the Dendriform operad, hence also to the
symmetric version of the Mould operad, one gets a pre-Lie product on
the free $\mld$ algebra on one generator and an injective map from the
free Pre-Lie algebra on one generator to the free $\mld$ algebra on
one generator. 

The pre-Lie product is defined by the formula
\begin{equation}
  f \sous g =\left(\frac{1}{u_1 u_{1\dots 2}}\circ_2 f\right) \circ_1
  g-\left(\frac{1}{u_{1\dots 2} u_2}\circ_2 g\right) \circ_1 f.
\end{equation}
or just as $g \succ f - f\prec g $. More explicitly, it is given by
\begin{multline}
  f(u_1,\dots,u_m)g(u_{m+1},\dots,u_{m+n})\frac{u_{m+1\dots
  m+n}}{u_{1\dots m+n}}\\- g(u_1,\dots,u_n)f(u_{n+1},\dots,u_{m+n})\frac{u_{1\dots
  n}}{u_{1\dots m+n}}.
\end{multline}

\begin{theorem}
  \label{sous_alternal}
  The pre-Lie product $\sous$ preserves alternality: if $f$ and $g$
  are alternal, then so is $f \sous g$.
\end{theorem}

\begin{proof}
  Let $f\in\mld(m)$ and $g\in\mld(n)$. Let us fix $j \in
  \{1,\dots,m+n-1\}$. One has to check that
  \begin{equation}
    \sum_{\sigma\in\sh(j,m+n-j)} (f\sous g)(u_{\sigma(1)},\dots,u_{\sigma(m+n)})=0.
  \end{equation}
  
  Let us compute the first half of this sum:
  \begin{equation}
     \sum_{\sigma\in\sh(j,m+n-j)} (g\succ f)(u_{\sigma(1)},\dots,u_{\sigma(m+n)}).
  \end{equation}
  This is
  \begin{multline}
    \frac{1}{u_{1\dots m+n}} \sum_{\sigma\in\sh(j,m+n-j)}
    f(u_{\sigma(1)},\dots,u_{\sigma(m)})\\
    g(u_{\sigma(m+1)},\dots,u_{\sigma(m+n)})(u_{\sigma(m+1)}+\dots+u_{\sigma(m+n)}). 
  \end{multline}
  Let us introduce the set $E'=\{\sigma(1),\dots,\sigma(m)\}$ and let
  $E''$ be its complement. Then, by standard properties of
  shuffles, one can rewrite the previous sum as
  \begin{equation}
    \frac{1}{u_{1\dots m+n}} \sum_{E'} \sum_{\sigma'}
    f(u_{\sigma'})
    \left(\sum_{k \in E''} u_{k} \right) \sum_{\sigma''} g(u_{\sigma''}),
  \end{equation}
  where $E'$ runs over the set of subsets of $\{1,\dots,m+n\}$ of
  cardinal $m$, $\sigma'$ is a shuffle of $E' \cap \{1,\dots ,j\}$ and
  $E'\cap \{j+1,\dots ,m+n\}$ and $\sigma''$ is a shuffle of $E'' \cap
  \{1,\dots ,j\}$ and $E''\cap \{j+1,\dots ,m+n\}$. Here we have used
  the shorthand notation introduced at the end of \S \ref{rappels}.

  Using the alternality of $f$ and $g$, this sum reduces to the terms
  where $E'$ is either included in $\{1,\dots ,j\}$ or in $\{j+1,\dots
  ,m+n\}$ and $E''$ is either included in $\{1,\dots ,j\}$ or in
  $\{j+1,\dots ,m+n\}$. Hence, the sum reduces to
  \begin{multline}
    \frac{1}{u_{1\dots m+n}} (f(u_1,\dots,u_j)g(u_{j+1},\dots
    ,u_{m+n})(u_{j+1}+\dots +u_{m+n})\\+f(u_{j+1},\dots
    ,u_{m+n})g(u_1,\dots,u_j)(u_1+\dots+u_j)),
  \end{multline}
  where the first term is present if $j=m$ and the second one if $j=n$.

  One can compute in the same way the other half of the sum:
  \begin{equation}
        \sum_{\sigma\in\sh(j,m+n-j)} (f\prec g)(u_{\sigma(1)},\dots,u_{\sigma(m+n)})
  \end{equation}
  and find exactly the same result. Hence the full sum vanishes as
  expected.
\end{proof}

\begin{corollary}
  The image of the free pre-Lie algebra on one generator is contained
  in the intersection of the free dendriform algebra with set of
  alternal elements in $\mld$.
\end{corollary}

It seems moreover that this inclusion could be an equality.

\begin{remark}\rm
  It follows also from the definition of $\sous$ given above that the
  Lie bracket $\mathrm{LIMU}$ associated to the associative product
  $\mathrm{MU}$ is also the bracket associated to the pre-Lie product
  $\sous$.
\end{remark}

\section{Suboperads in the category of sets}

The aim of this section is to describe two small suboperads of the
image of $\dend$ in $\mld$. The key point is that we work here in the
category of sets rather than in the category of vector spaces as
usual.

\subsection{Combinatorics of non-crossing trees and plants}

\begin{figure}
  \begin{center}
    \epsfig{file=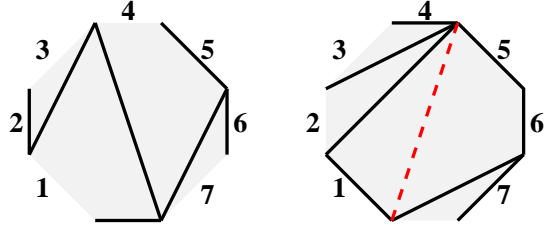,height=3cm} 
    \caption{a non-crossing tree and a non-crossing plant}
    \label{nct}
  \end{center}
\end{figure}

Let $n\geq 2$ be an integer. Consider the set of vertices of a regular
polygon with $n+1$ sides. One of these sides will be placed at the
bottom and called the \textit{base side}. The other sides are then
numbered from $1$ to $n$ from left to right. A \textit{diagonal} is a
line segment drawn between any two vertices of the regular polygon.
Two diagonals are \textit{crossing} if they are distinct and meet at
some point in the interior of the convex polygon.

An \textit{non-crossing plant} consists of two disjoint subsets of the
set of diagonals: the set of \textit{numerator diagonals} (pictured
dashed and red) and the set of \textit{denominator diagonals}
(pictured plain and black) with the following properties:
\begin{itemize}
\item any two diagonals in the union of these subsets are non-crossing,
\item the simplicial complex made by the denominator diagonals is
  connected and contains all vertices,
\item any numerator diagonal is contained in a closed cycle of
  denominator diagonals,
\item any closed cycle of denominator diagonals contains exactly one
  numerator diagonal.
\end{itemize}

In the sequel, a diagonal will always mean implicitly a denominator
diagonal, unless explicitly stated otherwise. Note that a side can
only be a denominator diagonal.

Let us call \textit{based non-crossing plant} a non-crossing plant
that includes the base side of the regular polygon.

If there is no numerator diagonal, then a non-crossing plant is a
\textit{non-crossing tree}, i.e. a maximal set of pairwise
non-crossing diagonals whose union is a connected and simply connected
simplicial complex (i.e. a tree).

Non-crossing trees and based non-crossing trees are well-known
combinatorial objects, see \cite{noy} for example. Non-crossing plants
seem not to have been considered before. Fig. \ref{nct} displays a
based non-crossing tree on the left and a non-crossing plant on the
right.

\begin{figure}
  \begin{center}
    \epsfig{file=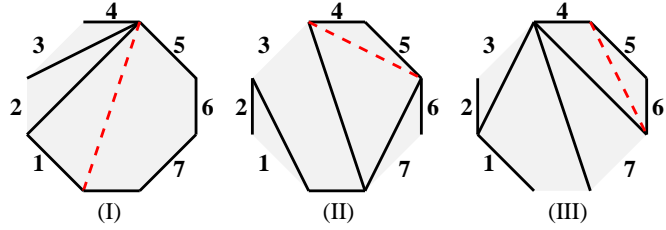,height=3cm} 
    \caption{the 3 sorts of non-crossing plants}
    \label{trio}
  \end{center}
\end{figure}

We need a precise recursive description of non-crossing plants. One
has to distinguish three sorts of them, as depicted in Fig.
\ref{trio}.

The first kind (I) is when the plant is based and the base side is
contained in a cycle of denominator diagonals, necessarily of length
at least $4$, as it must contain a numerator diagonal. For each other
diagonal in this cycle, one can consider diagonals that are in the
connected region (the one not containing the inner part of the cycle)
between this diagonal and the boundary of the regular polygon.  This
defines a based non-crossing plant. Conversely, one can pick any list
of length at least $3$ of based non-crossing plants and put them on
the sides of a closed cycle containing the base side and choose a
numerator diagonal in this cycle.

The second kind (II) is when the plant is based, but the base side is
not contained in a cycle of diagonals. Then there exists a unique side
which is not in the plant and which bounds the same region as the base
side. To the left and to the right of the square formed by this side
and the base side, one can define two non-crossing plants. This also
includes some degenerate cases, where one or both sides are empty. In
these cases, the square becomes a triangle or just the base side and
there is only one associated non-crossing plant (on the left or on the
right) or none.

The third and last kind (III) is when the plant is not based. Consider
the unique cycle (of length at least $3$) that would be created by
adding the base side. Just as in kind (I), one can define, by looking
at the outer regions bounded by this cycle, a list of based
non-crossing plants of length at least $2$.

Let us now translate this trichotomy in terms of generating series and
sketch the enumeration of non-crossing plants. We will use the
generating series
\begin{equation}
  P=\sum_{n\geq 1} p_n x^n=x+3 x^2+14 x^3+80 x^4+\dots,
\end{equation}
and
\begin{equation}
  Q=\sum_{n \geq 1} q_n x^n= x+2x^2+9x^3+51x^4+\dots,
\end{equation}
where $p_n$ (resp. $q_n$) is the number of non-crossing plants (resp.
based non-crossing plants) in the regular polygon with $n+1$ sides.

A non-crossing plant is either based or not. If it is not, it is of
kind $\mathrm{(III)}$ and can be described using a list (of length greater than
2) of based non-crossing plants. Hence one has
\begin{equation}
  P=Q+(Q^2/(1-Q))=Q/(1-Q).
\end{equation}

For a based non-crossing plant, either its base side is contained in a
cycle or not. In the first case, it is of kind $\mathrm{(I)}$ and can be
described using a list (of length $k$ at least 3) of based
non-crossing plants and the choice of an inner diagonal in a cycle of
length $k+1$. In the second case, it is of type $\mathrm{(II)}$ and can be
described by a pair of non-crossing plants or empty sides. Hence one
has
\begin{equation}
  Q=\sum_{k \geq 3} \frac{(k+1)(k-2)}{2}Q^k + x (1+P)^2.
\end{equation}

From these equations, one gets that $P$ satisfies the algebraic
relation:
\begin{equation}
  x-P+x P^2 +2xP +P^2+P^3=0.
\end{equation}
Therefore $P$ has a simple functional inverse:
\begin{equation}
  x=\frac{P-P^2-P^3}{(1+P)^2}.
\end{equation}
One can remark that this series appear as example (g) in \cite{loday_inverse}.

Let us now introduce the following notions, used in the next section.

A \textit{peeling point} of a non-crossing plant is a vertex (not in the
base side) such that the only incident diagonals are sides. In Figure
\ref{nct}, the non-crossing plant on the right has 2 peeling points between
sides 3 and 4 and between sides 5 and 6.
\begin{lemma}
  \label{peeling}
  There is always at least one peeling point in a non-crossing plant
  in the $n+1$ polygon, for $n\geq 2$.
\end{lemma}

\begin{proof}
  By induction on $n\geq 2$. This is true for all $3$ non-crossing
  plants $\gA$, $\dA$, $\mA$ in a triangle by inspection. Let us
  distinguish three cases, as before.

  (I) The base side belongs to a cycle of diagonals. 

  If there is something else than the cycle, there is a peeling point
  by induction in one of the sub-non-crossing plants bounding the
  cycle. This gives a peeling point in the whole non-crossing plant.

  If there is just a cycle, one can pick any vertex not in the base
  and not contained in the numerator diagonal. This vertex is a
  peeling point.

  (II) The base side is a diagonal, but does not belong to a cycle of
   diagonals.

   One can consider the left or right sub-non-crossing plant, which
   contains a peeling point by induction. This provides a peeling
   point in the full non-crossing plant.

  (III) The base side is not a diagonal.

  If there is something else than the would-be cycle (see the
  description of kind (III)), there is a peeling point by induction in
  one of the sub-non-crossing plants bounding this would-be cycle.
  This is also a peeling point in the full non-crossing plant.

  If there is just a would-be cycle, one can pick as peeling point any
  vertex not in the base.

\end{proof}

A \textit{border leaf} is a peeling point in a based non-crossing tree
that has only one incident side. In Figure \ref{nct}, there are 3
border leaves between sides 2 and 3, sides 4 and 5 and sides 6 and 7
in the non-crossing tree on the left.

\begin{lemma}
  \label{borderleaf}
  There is always at least one border leaf in a based non-crossing
  tree in the $n+1$ polygon, for $n\geq 2$.
\end{lemma}

\begin{proof}
  By induction on $n \geq 2$. This is clear if $n=2$, for the based
  non-crossing trees $\gA$ and $\dA$. Assume that $n\geq 3$.  In any
  based non-crossing tree, one can define a left and a right subtree,
  as the connected components of the tree minus its base side (every
  based non-crossing tree is of kind (II) as a non-crossing plant). At
  least one of them is not empty. It is enough to prove that there is
  a border leaf in one of them. One can therefore assume for instance
  that only the right subtree is not empty.

  Then consider the leftmost diagonal (other than the base side)
  emanating from the right vertex of the base side. Either there is a
  non-empty non-crossing tree to its right, hence a border leaf inside
  it by induction, or it has nothing to its right (it is a side). In
  this case, one can build a new non-crossing tree by shrinking the
  base side to a point and taking the side to its right as new base
  side. By induction, this new non-crossing tree has a border leaf.
  This implies that the initial tree has one too.
\end{proof}
 
\subsection{The operad of non-crossing plants}

Each diagonal is mapped to a linear function in variables $u_1,\dots,
u_n$ as follows: the base side is mapped to $u_{1\dots n}$. The
other sides are mapped to $u_1,\dots,u_n$ in the clockwise order. For
diagonals which are not sides, one considers the half plane not
containing the base side, with respect to this diagonal. This diagonal
is mapped to the sum of the values of the sides that are in this
half-plane.

To each non-crossing plant, one can then associate a rational function
in $\mld$ which is the product of the linear functions associated to
its numerator diagonals divided by the product of the linear functions
associated to its denominator diagonals. For instance, for the
non-crossing tree in the example of Fig. \ref{nct}, one gets
\begin{equation}
  1/ \left( u_{1\dots 7} u_2  u_{2\dots 3} u_{4\dots 7} u_5 u_6 u_{6\dots 7}  \right),
\end{equation}
and for the non-crossing plant on the right of the same figure,
\begin{equation}
 u_{1\dots 4}/\left( u_{1\dots 6} u_1 u_{2\dots 4}u_{3\dots 4} u_4 u_5 u_6 u_7  \right).
\end{equation}
The $3$ non-crossing plants in a triangle $\gA$,$\dA$ and $\mA$ are
mapped to $1/(u_1 u_{1\dots 2})$, $1/(u_{1\dots 2} u_2)$ and $1/(u_1
u_{2})$ in $\mld$.

This mapping from the set of non-crossing plants to $\mld$ is obviously
injective, as one can recover the non-crossing plant from the factorization
of its image. Therefore, we will from now on identify non-crossing plants
with their images in $\mld$.

\begin{figure}
  \begin{center}
    \epsfig{file=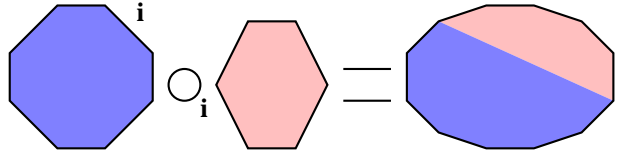,height=2cm} 
    \caption{grafting of non-crossing plants: deformation of polygons}
    \label{graft}
  \end{center}
\end{figure}

One can check, using the the definition of the composition in $\mld$,
that the set of non-crossing plants is closed under composition. Let
us give a combinatorial description of the composition of non-crossing
plants. Given two non-crossing plants $f$ and $g$ in some regular
polygons and a side $i$ of the polygon containing $f$, one has to
define a new non-crossing plant in the grafted polygon as in Fig
\ref{graft}. This is simply the union of $f$ and $g$, with some
modification along the grafting diagonal. If this diagonal is present
in both $f$ and $g$, then it is kept in $f\circ_i g$. If it is present
in exactly one of $f$ and $g$, then it is not kept in $f\circ_i g$. If
it is present in neither $f$ or $g$, then it becomes a numerator
diagonal in $f\circ_i g$.

\begin{theorem}
  The non-crossing plants form a suboperad $\NCP$ in the category of
  sets, which is contained in the image of $\dend$. This operad has
  the following presentation: three generators $\gA$, $\dA$ and $\mA$,
  subject only to the relations:
  \begin{align}
    \label{rela1}
    \gA \circ_2 \dA&= \dA \circ_1 \gA,\\
    \label{rela2}
    \mA \circ_1 \mA&= \mA \circ_2 \mA,\\
    \label{rela3}
    \gA \circ_1 \mA&=\gA \circ_2 \gA,\\
    \label{rela4}
    \dA \circ_2 \mA&=\dA\circ_1 \dA.
  \end{align}
\end{theorem}

\begin{proof}
  The inclusion in the image of $\dend$ will follow from the
  presentation by generators and relations, with generators $\gA$,
  $\dA$ and $\mA$ that are the images of the elements $\lft$, $\rgt$
  and $\lft+\rgt$ of $\dend$ in $\mld$.

  There remains to prove the presentation. Let us define an operad
  $\NCP'$ by the presentation above, with $\gA$, $\dA$ and $\mA$
  replaced by symbols $L$, $R$ and $M$. As the relations are satisfied
  in $\NCP$, there is a unique morphism of operads $\nabla$ from $\NCP'$
  to $\NCP$ sending the generators $L$, $R$ and $M$ to $\gA$, $\dA$ and
  $\mA$. Let us prove by induction that there is an inverse $\Delta$
  to $\nabla$.

  Note that, whenever this makes sense, $\Delta$ is of course a
  morphism of operads. The existence of $\Delta$ is clear for $n=2$.
  Assume that $n\geq 3$ and let $T$ be in $\NCP(n)$. By Lemma
  \ref{peeling}, there is a peeling point in the non-crossing plant
  $T$. Let $i$ be the index of the side of the polygon which is left
  to this leaf.  Then $T$ can be written $S \circ_i \delta$ where
  $\delta$ is either $\dA$, $\gA$ or $\mA$ and $S \in \NCP(n-1)$.

  Let us define
  \begin{equation}
    \Delta(T)=\Delta(S)\circ_i \Delta(\delta).
  \end{equation}

  One has to prove that this definition does not depend on the choice
  of the peeling point. Let us assume that there is another peeling
  point.  Without further restriction, one can assume that it is at
  the right of side $j$ with $i < j$. Thus $T$ can also be written
  $S'\circ_j \delta'$ where $\delta'$ is either $\gA$, $\dA$ or $\mA$
  and $S' \in \NCP(n-1)$. One then has to distinguish two cases.

  \textbf{Far case}

  If $i+1 < j$, then there is still a peeling point in
  $S'$ at the right of edge $i$, hence there exists $S''\in \NCP(n-2)$
  such that $T$ can be written as
  \begin{equation}
    S'\circ_j \delta' =(S'' \circ_i \delta) \circ_j \delta'
=(S''\circ_{j-1} \delta') \circ_i \delta=S \circ_i \delta,
  \end{equation}
  where the second equality is an axiom of operads. This implies that
  both choices of peeling point in $T$ leads to the same value for
  $\Delta(T)$:
  \begin{multline}
    \Delta(S')\circ_j \Delta(\delta')=(\Delta(S'') \circ_i \Delta(\delta)) \circ_j \Delta(\delta')\\=(\Delta(S'')\circ_{j-1} \Delta(\delta'))
    \circ_i \Delta(\delta)=\Delta(S) \circ_i \Delta(\delta).
  \end{multline}

  \textbf{Near cases}

  If $i+1=j$, then one can distinguish $4$ cases. Either all three
  sides $i,i+1,i+2$ are diagonals, or just two of them are. The other
  possibilities are excluded by the second condition in the definition
  of non-crossing plants.

  Let us consider the first case. Necessarily $T$ can be written, for
  some $S''\in \NCP(n-2)$, as
  \begin{equation}
    (S''\circ_i \mA) \circ_{i+1} \mA=S''\circ_i (\mA \circ_2
    \mA)=S''\circ_i (\mA \circ_1 \mA)=(S''\circ_i \mA) \circ_i \mA. 
  \end{equation}
  This implies that both choices of peeling point give the same value
  for $\Delta(T)$:
  \begin{equation}
    \Delta(S''\circ_i \mA) \circ_{i+1} M=(\Delta(S'')\circ_i M) \circ_{i+1} M=(\Delta(S'')\circ_i M) \circ_i M
    =\Delta(S''\circ_i \mA) \circ_i M, 
  \end{equation}
  where the middle equality follows from relation (\ref{rela2}) for
  the generators $L$, $R$ and $M$.

  The three other cases are similar to this one, each one of them
  using one of the relations (\ref{rela1}), (\ref{rela3}) and
  (\ref{rela4}) for the generators $L$, $R$ and $M$.

  Hence $\Delta$ is well-defined. Then, one has for all $T$ in $\NCP(n)$,
  \begin{equation}
    \nabla(\Delta(T))=\nabla(\Delta(S) \circ_i \Delta(\delta))=S \circ_i \delta=T,
  \end{equation}
  by induction hypothesis and because $\nabla$ is a morphism of
  operad. Let $x$ be in $\NCP'(n)$. Then $x$ can written $y \circ_i d$
  for some $y$ in $\NCP'(n-1)$ and $d \in \{L,R,M\}$. Then
  \begin{equation}
    \Delta(\nabla(x))=\Delta(\nabla(y) \circ_i \nabla(d))
    =\Delta(\nabla(y)) \circ_i \Delta(\nabla(d))=y \circ_i d=x.
  \end{equation}
  Here for the computation of $\Delta$ we choose the peeling point
  corresponding to $\nabla(d)$. We have proved that $\Delta$ is the
  inverse of $\nabla$ up to order $n$. This concludes the induction step.

\end{proof}

\begin{remark}\rm
  By adding the opposite of each non-crossing plant, one can get an
  anticyclic operad in the category of sets. The cyclic action is just
  given by the rotation of the regular polygon, up to sign.
\end{remark}

\subsection{The operad of based non-crossing trees}

It is not hard to see, using the combinatorial description of the
composition given above, that based non-crossing trees are closed under
composition.

\begin{proposition}
  The suboperad of $\NCP$ generated by $\gA$ and $\dA$ is exactly the
  suboperad of based non-crossing-trees. These generators are only
  subject to the relation {\rm (\ref{rela1})}.
\end{proposition}

\begin{proof}
  As $\gA$ and $\dA$ are based non-crossing trees, the operad they
  generate is contained in the suboperad of non-crossing trees.

  To prove the reverse inclusion, one proceeds by induction. By Lemma
  \ref{borderleaf}, there is a border leaf in any non-crossing tree
  $T$. Let $i$ be the index of the side of the polygon which is left
  to this leaf of $T$. Then $T$ can be written $S \circ_i \delta$
  where $\delta$ is either $\dA$ or $\gA$ and $S$ is a smaller based
  non-crossing tree. This implies the inclusion in the suboperad
  generated by $\dA$ and $\gA$.

  The presentation is a consequence of that of the bigger operad
  $\NCP$.
\end{proof}

From results of Loday in \cite{atree}, one can deduce that
\begin{proposition}
  For any non-crossing plant $T$ of degree $n$, the inverse image of
  $T$ in $\dend$ is a sum without multiplicities in $\dend(n)$ and can
  therefore be considered as a subset of $\bY(n)$.
\end{proposition}

We would like to draw the attention on the following conjecture, which
has been checked in low degrees. Recall that the Tamari poset
\cite{tamari} is a partial order on the set $\bY(n)$ which indexes a
basis of $\dend(n)$.

\begin{conjecture}
  For any non-crossing tree $T$ of degree $n$ (not necessarily based)
  in $\mld$, the inverse image of $T$ in $\dend$ is
  \begin{equation}
    \sum_{t \in I} t,
  \end{equation}
  where $I$ is some interval in the Tamari poset $\bY(n)$.
\end{conjecture}

\section{Other structures}

Let us consider some other operations on the free $\mld$ algebra on
one generator.

\subsection{Over and Under operations}

Loday and Ronco have introduced in \cite{lodayronco} two other
associative products on the free dendriform algebra on one generator,
called \textit{Over} and \textit{Under} and denoted by $/$ and
$\backslash$. They are usually defined as simple combinatorial
operations on planar binary trees, but can be restated using the
Dendriform operad as follows:
\begin{equation}
  f / g = (g\circ_1 \lft)\circ_1 f
\end{equation}
and
\begin{equation}
  f \backslash g = (f \circ_m \rgt )\circ_{m+1} g,
\end{equation}
where $f$ is assumed to be of degree $m$. One can use this to extend
these operations to the free $\mld$ algebra on one generator.
Explicitly, restated inside the $\mld$ operad, these products are
given by
\begin{equation}
  (f / g) (u_1,\dots,u_{m+n})=f(u_1,\dots,u_n)g(u_{1\dots n+1},u_{n+2},\dots,u_{m+n})
\end{equation}
and
\begin{equation}
  (f \backslash g) (u_1,\dots,u_{m+n})=f(u_1,\dots,u_{n-1},u_{n\dots m+n})g(u_{n+1},\dots,u_{m+n}).
\end{equation}

\subsection{Structures associated to the operad structure}

The fact that $\mld$ is an operad implies that one can define other
operations on the free $\mld$ algebra on one generator, namely a
pre-Lie product $\circ$ (not to be confused with the one introduced
before and denoted by $\sous$) and the associated Lie bracket and
group law.

The pre-Lie product $\circ$ is defined for $f\in\mld(m)$ and
$g\in\mld(n)$ by
\begin{equation}
  f \circ g= \sum_{i=1}^{m} f\circ_i g.
\end{equation}
More explicitly, $f\circ g$ is given by
\begin{equation}
  \sum_{i=1}^{m} f(u_1,\dots,u_{i-1},u_{i\dots
  i+n-1},u_{i+n},\dots,u_{m+n-1})g(u_i,\dots ,u_{i+n-1}) u_{i\dots
  i+n-1}.
\end{equation}

This construction is clearly functorial from the category of operads
to the category of pre-Lie algebras. Note that the product $f \circ g$
is in $\mld(m+n-1)$.

\begin{theorem}
  The pre-Lie product $\circ$ preserves alternality, that is $f
  \circ g$ is alternal as soon as $f$ and $g$ are.
\end{theorem}

\begin{proof}
  Let $f$ be in $\mld(m)$ and $g\in\mld(n)$. Let us fix $j \in
  \{1,\dots,m+n-2\}$. One has to check that
  \begin{equation}
    \sum_{\sigma\in\sh(j,m+n-1-j)}(f\circ g)(u_{\sigma(1)},\dots,u_{\sigma(m+n-1)})=0.
  \end{equation}
  The sum to be computed is
  \begin{multline}
    \sum_{i=1}^{m} \sum_{\sigma\in\sh(j,m+n-1-j)}
    f(u_{\sigma(1)},\dots,u_{\sigma(i)}+\dots +
    u_{\sigma(i+n-1)},\dots,u_{\sigma(m+n-1)})\\g(u_{\sigma(i)},\dots
    ,u_{\sigma(i+n-1)}) ( u_{\sigma(i)}+\dots+ u_{\sigma(i+n-1)}).
  \end{multline}
  Let us introduce the sets $E'=\{\sigma(1),\dots,\sigma(i-1)\}$ of
  cardinal $i-1$ and $E''=\{\sigma(i+n),\dots,\sigma(m+n-1)\}$ of
  cardinal $m-i$. They
  have the following properties:
  \begin{itemize}
  \item $E' \cap \{1,\dots,j\}$ is an initial subset of $\{1,\dots,j\}$,
  \item $E' \cap \{j+1,\dots,m+n-1\}$ is an initial subset of
    $\{j+1,\dots,m+n-1\}$,
  \item $E'' \cap \{1,\dots,j\}$ is a final subset of $\{1,\dots,j\}$,
  \item $E'' \cap \{j+1,\dots,m+n-1\}$ is a final subset of $\{j+1,\dots,m+n-1\}$.
  \end{itemize}
  Let $E=\{\sigma(i),\dots,\sigma(i+n-1)\}$ be the complement of
  $E'\cup E''$. It follows from the conditions above that $E \cap
  \{1,\dots,j\}$ and $E \cap \{j+1,\dots,m+n-1\}$ are sub-intervals.

  Using standard properties of the set of shuffles, and the shorthand
  notation $u_\sigma$ introduced at the end of \S \ref{rappels}, one
  can then rewrite the previous sum as
  \begin{equation}
    \sum_{i=1}^{m} \sum_{E',E''} \sum_{\sigma',\sigma''} f(u_{\sigma'},\sum_{k\in E}
    u_k,u_{\sigma''} )  \sum_{\nu} g(u_{\mu}) \left(\sum_{k \in E} u_k\right),
  \end{equation}
  where $E'$ of cardinal $i-1$ and $E''$ of cardinal $m-i$ are subsets
  with the above properties, $\sigma'$ is a shuffle of $E'\cap
  \{1,\dots,j\}$ and $E' \cap \{j+1,\dots,m+n-1\} $, $\sigma''$ is a
  shuffle of $E'' \cap \{1,\dots,j\} $ and $E'' \cap
  \{j+1,\dots,m+n-1\} $ and $\nu$ is a shuffle of $E \cap
  \{1,\dots,j\} $ and $E \cap \{j+1,\dots,m+n-1\}$.

  Using the alternality of $g$, one can see that the sum reduces to
  the cases when $E \subset \{1,\dots,j\} $ or $E \subset
  \{j+1,\dots,m+n-1\}$. Let us show that each of these two terms
  vanishes. As the proof is similar, we treat only the case when $E
  \subset \{1,\dots,j\}$. In this case, there exists $k$ such that
  $E=\{k,\dots,k+n-1\}$. The corresponding term is
  \begin{equation}
    \sum_{i=1}^{m} \sum_{E',E''} \sum_{\sigma',\sigma''} f(u_{\sigma'}, u_{k\dots k+n-1},u_{\sigma''} ) g(u_{k},\dots,u_{k+n-1}) u_{k\dots k+n-1},
  \end{equation}
  where $E'$ and $E''$ runs over the appropriate sets.

  Once again by the usual properties of shuffles, this can be
  rewritten as
  \begin{equation}
    \sum_{k} g(u_k, \dots,u_{k+n-1}) (u_{k\dots k+n-1}) \sum_{\mu} f(u_{\mu}),
  \end{equation}
  where $1 \leq k \leq k+n-1 \leq j $ and $\mu$ is a shuffle of
  $\{1,\dots,{k-1},{(k \dots k+n-1)},{k+n},\dots,{j} \}$ and
  $\{j+1,\dots,m+n-1\}$, with the abuse of notation made in
  considering $(k \dots k+n-1)$ as an index. This sum is zero because
  $f$ is alternal.
\end{proof}

\begin{proposition}
  \label{circ_vegetal}
  If $f$ and $g$ are vegetal, then so is $f\circ g$.
\end{proposition}
\begin{proof}
  One has to compute
  \begin{multline}
 \sum_{\sigma\in\sym_{m+n-1}}
\sum_{i=1}^{m}     f(t\,
u_{\sigma(1)},\dots,t(u_{\sigma(i)}+\dots+u_{\sigma(i+n-1)}),\dots,t\,u_{\sigma(m+n-1)})\\
{u_1 \dots u_{m+n-1}} g(t\,u_{\sigma(i)},\dots,t\,u_{\sigma(i+n-1)})
t\, (u_{\sigma(i)}+\dots+u_{\sigma(i+n-1)}).
  \end{multline}
  Let us introduce the set $E=\{\sigma(i),\dots,\sigma(i+n-1)\}$. One
  can then rewrite the previous sum as
  \begin{multline}
      {u_1 \dots u_{m+n-1}\,t}\sum_E (\sum_{j\in E} u_j) 
     \sum_{i=1}^{m} \sum_{\sigma'} f(t\,
u_{\sigma'(1)},\dots,t(\sum_{j\in E} u_j),\dots,t\,u_{\sigma'(m+n-1)}) \\\sum_{\sigma''} g(t\,u_{\sigma''(i)},\dots,t\,u_{\sigma''(i+n-1)}),
  \end{multline}
  where $E$ runs over the set of subsets of cardinal $n$ of
  $\{1,\dots,m+n-1\}$, $\sigma'$ is a bijection from
  $\{1,\dots,i-1\}\sqcup \{i+n,\dots,m+n-1\} $ to
  the complement of $E$ and $\sigma''$ is a bijection from
  $\{i,\dots,i+n-1\}$ to $E$.

  Using first the vegetal property of $g$, one gets
  \begin{multline}
    {u_1 \dots u_{m+n-1} t}\sum_E \frac{\sum_{j\in E} u_j }{\prod_{j
        \in E} u_{j}} \sum_{i=1}^{m} \sum_{\sigma'} f(t\,
    u_{\sigma'(1)},\dots,t\left(\sum_{j\in E} u_j\right),\dots,t\,u_{\sigma'(m)})
    \\n! g(t,\dots,t).
  \end{multline}
  Then, this can be rewritten as
  \begin{equation}
    {n! t} 
    g(t,\dots,t) \sum_E \left(\sum_{j\in E} u_j\right) 
    \sum_{\theta} f(t\,
    \theta(1),\dots,t\, \theta(m)) \prod_{j \not \in E} u_{j}.
  \end{equation}
  where $\theta$ is a bijection from $\{1,\dots,m\}$ to $\{u_j\}_{j
    \not\in E} \sqcup \{\sum_{j\in E} u_j\}$. By using the vegetal
  property of $f$, this becomes
  \begin{equation}
    {n!}\binom{m+n-1}{n} m! f(t,\dots,t) g(t,\dots,t) t.
  \end{equation}
  Using once again the vegetal property of $f$ in the special case
  $u_1=n$ and $u_2=\dots=u_m=1$, one gets
  \begin{equation}
    (m+n-1)! \sum_{i=1}^{m} f(t,\dots,n t,\dots, t)\, g(t,\dots,t) \,nt,
    \end{equation}
  which is $(m+n-1)! (f\circ g)(t,\dots,t)$. Hence $f\circ g$ is vegetal.
\end{proof}

\subsection{Forgetful morphism to formal vector fields}

The aim of this section is to define a map $\forg$ of pre-Lie algebras
from moulds satisfying appropriate conditions to formal power series
in one variable $x$, or rather to formal vector fields.

\smallskip

Let us consider here only moulds $f$ such that $f_n$ is homogeneous of
weight $-n$ and has only poles at some $u_{i\dots j}$ with arbitrary
multiplicity (nice poles). From Remarks \ref{homogen} and
\ref{nicepole}, this subspace is a anticyclic suboperad, hence it is
closed for $\circ$ by functoriality. Let us consider its intersection
with the subspace of vegetal moulds, which is also closed for $\circ$
by Prop. \ref{circ_vegetal}.  Let us note that this intersection
contains the image of $\dend$ by Prop.  \ref{dend_vegetal} and Prop.
\ref{fraction_arbre}.

Let us recall that the usual pre-Lie product, also denoted by $\circ$,
on vector fields is given by
\begin{equation}
  F(x)\partial_x \circ G(x)\partial_x = (\partial_x F(x)) G(x) \partial_x.
\end{equation}

\begin{theorem}
  The substitution $u_i \mapsto 1/x$ induces a morphism $\forg$ of
  pre-Lie algebras $f \mapsto f(x^{-1},\dots,x^{-1})\partial_x$ from
  $\mld$ (restricted as above to homogeneous vegetal moulds with nice
  poles) with the pre-Lie product $\circ$ to the pre-Lie algebra of
  vector fields in the variable $x$ with formal power series in $x$ as
  coefficients.
\end{theorem}

\begin{proof}
  Let $f\in \mld(m)$ and $g\in\mld(n)$ satisfying the additional
  conditions stated before. One has to prove that
  \begin{multline}
    \partial_x f(x^{-1},\dots,x^{-1})g(x^{-1},\dots,x^{-1})\\
    =\sum_{i=1}^{m} n x^{-1}  f(x^{-1},\dots,n
    x^{-1},\dots,x^{-1})g(x^{-1},\dots,x^{-1}),
  \end{multline}
  where we have used the assumed shape of the poles of $f$ and $g$ to
  ensure that the substitution makes sense. It is therefore enough to
  prove that
  \begin{equation}
    x \partial_x f(x^{-1},\dots,x^{-1})
    =\sum_{i=1}^{m} n f(x^{-1},\dots,n
    x^{-1},\dots,x^{-1}).
  \end{equation}
  From the homogeneity of $f$, one has to prove that
  \begin{equation}
    m f(x^{-1},\dots,x^{-1})=\sum_{i=1}^{m} n f(x^{-1},\dots,n
    x^{-1},\dots,x^{-1}),
  \end{equation}
  which is a special case of the vegetal property of $f$, with
  $t=x^{-1}$, $u_1=n$ and $u_2=\dots=u_m=1$.
\end{proof}

As the group law associated to the usual pre-Lie product on formal
power series is the classical composition of power series, one can see
the group structure on moulds corresponding to $\circ$ as some kind of
generalized composition.

\begin{remark}\rm
  From (\ref{def_MU}), it is quite obvious that the associative
  product $\mathrm{MU}$ is mapped by $\forg$ to the usual commutative
  product of formal power series.
\end{remark}

\subsection{Derivation}

Let us introduce a map $\partial$ on moulds, which decreases the
degree by $1$.

For a mould $f \in \mld(m)$, $\partial f$ is the element of
$\mld(m-1)$ defined by
\begin{equation}
  \partial f (u_1,\dots,u_{m-1})= \sum_{j=1}^m \res_{t=0} f(u_1,\dots,u_{j-1},t,u_j,\dots,u_{m-1}),
\end{equation}
where $\res$ is the residue.

The main motivation for this map is the following property.
\begin{proposition}
  The map $\partial$ is sent by the forgetful map $\forg$ to the
  partial derivative with respect to $x$, i.e. for any $f \in \mld(m)$
  which is homogeneous, vegetal and has nice poles, one has
  $\forg(\partial f)=\partial_x \forg(f)$.
\end{proposition}

\begin{proof}
  Let $f\in\mld(m)$. By homogeneity,
  $\forg(f)=f(x^{-1},\dots,x^{-1})=x^m f(1,\dots,1)$. Hence
  $\partial_x \forg(f)= m x^{m-1} f(1,\dots,1)$.

  On the other hand, $\forg(\partial f)$ is
  \begin{equation}
    \res_{t=0} \sum_{j=1}^m f(x^{-1},\dots,t,\dots,x^{-1}),
  \end{equation}
  where $t$ is in the $j$th position. By the vegetal property of $f$,
  this is
  \begin{equation}
    \res_{t=0} m \frac{x^{m-1}}{t} f(1,\dots,1).
  \end{equation}
  This proves the expected equality.
\end{proof}

\begin{remark}\rm
  If $f$ is alternal and of degree at least $2$, then $\partial f=0$.
  This is obvious once the definition of $\partial$ is rewritten as
  the residue of a sum over shuffles of $t$ with
  $\{u_1,\dots,u_{m-1}\}$.
\end{remark}

\begin{proposition}
  \label{partial_deri}
  The map $\partial$ is a derivation for the products
  $\prec,\succ,\sous$, $\mathrm{MU}$ and $\mathrm{LIMU}$. It is also a
  derivation for $\circ$, under the restriction that functions have
  nice poles.
\end{proposition}

\begin{proof}
  It is enough to prove this for $\succ$ and $\circ$. The case of
  $\prec$ is similar to the case of $\succ$ and the other cases can be
  deduced from these ones.

  Let us consider the case of $\succ$. Let $f\in\mld(m)$ and
  $g\in\mld(n)$. One has to compute
  \begin{multline}
    \sum_{j=1}^{m} \res_{t=0} \left(f(u_1,\dots,t,u_j,\dots, u_{m-1})g(u_{m},\dots,u_{m+n-1})\frac{u_{m\dots
          m+n-1}}{u_{1\dots m+n-1}+t} \right) + \\\sum_{j=m+1}^{m+n} \res_{t=0} \left( f(u_1,\dots, u_m)g(u_{m+1},\dots,t,u_j,\dots,u_{m+n-1})\frac{u_{m+1\dots
          m+n-1}+t}{u_{1\dots m+n-1}+t}.
    \right).
  \end{multline}

  By the properties of the residue, this becomes
  \begin{multline}
    \sum_{j=1}^{m} \res_{t=0} \left(f(u_1,\dots,t,u_j,\dots,
      u_{m-1})\right) g(u_{m},\dots,u_{m+n-1})\frac{u_{m\dots
          m+n-1}}{u_{1\dots m+n-1}}  + \\\sum_{j=m+1}^{m+n}
      f(u_1,\dots, u_m) \res_{t=0}
      \left(g(u_{m+1},\dots,t,u_j,\dots,u_{m+n-1})\right) \frac{u_{m+1\dots
          m+n-1}}{u_{1\dots m+n-1}}.    
  \end{multline}
  This is $\partial f \succ g + f \succ \partial g$, which proves that
  $\partial$ is a derivation of $\succ$.

  Let us now consider the case of $\circ$. One has to compute
  \begin{multline}
    \sum_{i=1}^{m} \sum_{j=1}^{i-1} \res_{t=0} (
      f(u_1,\dots,t,u_j,\dots, u_{i-2},u_{i-1\dots
  i+n-2},u_{i+n-1},\dots,u_{m+n-2})\\g(u_{i-1},\dots ,u_{i+n-2}) u_{i-1\dots
  i+n-2}
)+ \\ \sum_{i=1}^{m} \sum_{j=i}^{i+n-1}\res_{t=0} ( f(u_1,\dots,u_{i-1},u_{i\dots
  i+n-2+t},u_{i+n-1},\dots,u_{m+n-2})\\g(u_i,\dots,t,u_j,\dots ,u_{i+n-2}) (u_{i\dots
  i+n-2}+t)) + \\ \sum_{i=1}^{m} \sum_{j=i+n}^{m+n-1} \res_{t=0} ( f(u_1,\dots,u_{i-1},u_{i\dots
  i+n-1},u_{i+n},\dots,t,u_j,\dots,u_{m+n-2})\\g(u_i,\dots ,u_{i+n-1}) u_{i\dots
  i+n-1}
).
  \end{multline}
  By the properties of residues, and the assumption that $f$ and $g$
  have nice poles, this becomes
  \begin{multline}
    \sum_{i=1}^{m} \sum_{j=1}^{i-1} \res_{t=0} \left(
      f(u_1,\dots,t,u_j,\dots, u_{i-2},u_{i-1\dots
  i+n-2},u_{i+n-1},\dots,u_{m+n-2})
\right) \\g(u_{i-1},\dots ,u_{i+n-2}) u_{i-1\dots
  i+n-2} + \\ \sum_{i=1}^{m} \sum_{j=i+n}^{m+n-1} \res_{t=0} \left( f(u_1,\dots,u_{i-1},u_{i\dots
  i+n-1},u_{i+n},\dots,t,u_j,\dots,u_{m+n-2})
\right) \\g(u_i,\dots ,u_{i+n-1}) u_{i\dots
  i+n-1}+ \\ \sum_{i=1}^{m} \sum_{j=i}^{i+n-1} f(u_1,\dots,u_{i-1},u_{i\dots
  i+n-2},u_{i+n-1},\dots,u_{m+n-2}) \\ \res_{t=0} \left( g(u_i,\dots,t,u_j,\dots ,u_{i+n-2})
\right) u_{i\dots
  i+n-2}.
  \end{multline}
  The first two terms give $\partial f \circ g$ and the third one
  gives $f \circ \partial g$.
\end{proof}

As a corollary of Prop. \ref{partial_deri}, the map $\partial$
preserves the image by $\psi$ of the free Dendriform algebra on one
generator. One can be more precise: the action of $\partial$ is by
vertex-removal, in the following sense. From the description of
$\psi(T)$ for a planar binary tree $T$ in Prop. \ref{fraction_arbre},
one can see that taking the residue with respect to one of the
variables and then renumbering correspond to the removal of a top
vertex in $T$. Hence $\partial \psi(T)$ is the sum over all top
vertices of $T$ of the image by $\psi$ of some smaller binary tree.

\subsection{The products $\mathrm{ARIT}$ and $\mathrm{ARI}$}

Let us define two other bilinear products on the free $\mld$ algebra
on one generator, denoted by $\mathrm{ARIT}$ and $\mathrm{ARI}$:
\begin{equation}
  \mathrm{ARIT}(f,g)=f \circ (g / Y) - f \circ (Y \backslash g)
\end{equation}
and
\begin{equation}
  \mathrm{ARI}(f,g)=\mathrm{ARIT}(f,g)-\mathrm{ARIT}(g,f)+\mathrm{LIMU}(f,g).
\end{equation}

One can check, by writing the explicit expression for these products,
that they do indeed reproduce the $\mathrm{ARIT}$ and $\mathrm{ARI}$
maps introduced by Ecalle in his study of moulds. In particular, it is
known that $\mathrm{ARI}$ is a Lie bracket that preserves alternality.
Note that we have defined as $\mathrm{ARIT}(f,g)$ what Ecalle denotes
by $\mathrm{ARIT}(g)f$.

\begin{lemma}
  \label{deri_ari}
  There holds $ \partial (f /Y) = (\partial f) / Y$ and
  $\partial (Y \backslash f) =  Y \backslash (\partial f)$.
\end{lemma}
\begin{proof}
  Let us consider only the first case, the other one being
  similar. Let $f\in\mld(n)$. As 
  \begin{equation}
    f / Y=\lft \circ_1 f=\frac{1}{u_{1\dots n+1}}f(u_1,\dots,u_n),
  \end{equation}
  one has to compute
  \begin{equation}
    \sum_{i=1}^{n} \res_{t=0} \frac{1}{u_{1\dots n}+t}
    f(u_1,\dots,u_{i-1},t,u_{i},\dots,u_{n-1})+\res_{t=0} \frac{1}{u_{1\dots n}+t}
    f(u_1,\dots,u_{n}).
  \end{equation}
  The second term vanishes, and what remains is
  \begin{equation}
     \frac{1}{u_{1 \dots n}}  \sum_{i=1}^{n} \res_{t=0} 
    f(u_1,\dots,u_{i-1},t,u_{i},\dots,u_{n-1}),
  \end{equation}
  which is exactly $ (\partial f) / Y $. 
\end{proof}

\begin{corollary}
  The map $\partial$ is a derivation for $\mathrm{ARIT}$ and
  $\mathrm{ARI}$.
\end{corollary}
\begin{proof}
  As $\partial$ is a derivation of $\circ$ by Prop. \ref{partial_deri}, one has
  \begin{equation}
    \partial(\mathrm{ARIT}(f,g))=\partial(f) \circ (g / Y)+f \circ \partial(g /
    Y) - \partial(f) \circ
    (Y \backslash g)- f \circ
    \partial(Y \backslash g).
  \end{equation}
  One can then conclude for $\mathrm{ARIT}$ using Lemma \ref{deri_ari}. The proof
  for $\mathrm{ARI}$ follows from this and from Prop. \ref{partial_deri}.
\end{proof}

Let us now state some properties of the $\mathrm{ARI}$ and
$\mathrm{ARIT}$ products.

\begin{proposition}
  The free dendriform algebra in $\mld$ is closed under
  $\mathrm{ARIT}$ and $\mathrm{ARI}$. The products $\mathrm{ARIT}$ and
  $\mathrm{ARI}$ preserves the vegetal property.
\end{proposition}

\begin{proof}
  For the first statement, it is enough to look at the definition of
  $\mathrm{ARIT}$. One already knows that the Over and Under operations are
  defined on the dendriform subspace. The product $\circ$ has the same
  property, because it is a functorial construction on operads.

  Let us prove the second statement. It is enough to prove this for
  $\mathrm{ARIT}$, by Prop. \ref{mu_vegetal}. As this property is already known
  for $\circ$ by Prop. \ref{circ_vegetal}, it is enough to prove that
  $ f /Y $ and $Y \backslash f$ are vegetal if $f$ is so. Let us
  consider only the first case, as the other one is just the same.

  Let $f$ be in $\mld(n)$. One has to compute
  \begin{equation}
    \sum_{\sigma \in \sym_{n+1}} f(u_{\sigma(1)},\dots,u_{\sigma(n)})\frac{1}{u_{1\dots n+1}\,t}.
  \end{equation}
  One can separate the sum according to the value of $\sigma(n+1)$,
  obtaining
  \begin{equation}
    \frac{1}{u_{1\dots n+1}\,t} \sum_{i=1}^{n+1} \sum_{\sigma'} f(u_{\sigma'(1)},\dots,u_{\sigma'(n)}),
  \end{equation}
  where $\sigma'$ runs over the bijections from $\{1,\dots,n\}$ to
  $\{1,\dots,n+1\}\setminus \{i\}$. Using then the vegetal property of
  $f$, this becomes
  \begin{equation}
       \frac{1}{u_{1\dots n+1}\,t} \sum_{i=1}^{n+1} \frac{u_i}{u_1\dots
         u_{n+1}} n! f(t,\dots,t).
  \end{equation}
  This gives
  \begin{equation}
    \frac{f(t,\dots,t) (n+1)!}{ u_1\dots
         u_{n+1}\, (n+1) \, t },
  \end{equation}
  which proves the expected vegetal property.
\end{proof}

\section{Examples of moulds}

Let us describe the image in $\mld$ of some special and nice elements
of $\dend$.

\smallskip

Let $AS$ be the mould defined by $AS_n=1/(u_1\dots u_n)$. The components
of this mould provide the basis of the Associative suboperad. Hence it
is in the image of the Dendriform operad in the Mould operad. On the
other hand, it is known that the basis of the Associative suboperad of
the Dendriform operad is given by the sum of all planar binary trees.
Hence, one has
\begin{equation}
  AS_n=\psi\left(\sum_{t\in\bY(n)} T\right)=\frac{1}{u_1 \dots u_n}.
\end{equation}
On can note that the image by $\forg$ of the mould $AS$ is $x/(1-x)$.

\smallskip

Let us say that a planar binary tree is of type (p,q) if its left
subtree has $p+1$ leaves and its right subtree has $q+1$ leaves. The
sum over binary trees of type (p,q) is 
\begin{equation}
  \frac{u_p}{u_1 \dots u_n
  (u_{1\dots n}) }.
\end{equation}
This is an easy consequence of the previous result, using for instance
the Over and Under products.

Let $TY$ be the mould defined by 
\begin{equation}
  TY_n=\frac{\sum_{i=1}^n t^{i-1} u_i}{u_1 \dots u_n u_{1\dots n}},
\end{equation}
with a parameter $t$. By the preceding discussion, the mould $TY$ is
also in the image of the Dendriform operad. The image of $TY$ by
$\forg$ is
\begin{equation}
  \frac{1}{1-t}\log\left(\frac{1-tx}{1-x}\right).
\end{equation}

Another interesting mould has the following components:
\begin{equation}
  \frac{\sum_{i=1}^{n} i \, u_i}{u_1 \dots u_n u_{1\dots n}}.
\end{equation}
By the same argument as above, this mould belongs to the image of the
Dendriform operad. This mould should be related to the series indexed
by planar binary trees considered in \cite{int_tamari}. Its image by
$\forg$ is
\begin{equation}
  \frac{x(2-x)}{2 (1-x)^2}.
\end{equation}

\smallskip

One can also compute the image of the Connes-Moscovici series. Let us
first recall its definition. In the free pre-Lie algebra on one
generator, where the product is denoted by $\sous$, let $CM_1$ be the
generator and let
\begin{equation}
  CM_n= CM_{n-1} \sous CM_1.
\end{equation}

One can consider these objects as elements of the free dendriform
algebra on one generator, endowed with the pre-Lie product $\sous$. It
follows from Prop. \ref{sous_alternal} that these elements are
alternal.

\begin{proposition}
  One has 
  \begin{equation}
    \psi(CM_n)=\frac{1}{u_1 \dots u_n u_{1\dots n}}\sum_{k=1}^n
    (-1)^{n+k}\binom{n}{k} u_k.
  \end{equation}
  The image of $\psi(CM)$ by $\forg$ is $x$.
\end{proposition}

\begin{proof}
  The proof is by induction. By definition of the pre-Lie operation
  $\sous$ in the Dendriform operad, one has
  \begin{equation}
    CM_{n+1}= \lft \circ_2 CM_{n} - \rgt \circ_1 CM_{n}.
  \end{equation}
  Hence, in $\mld$, one gets
  \begin{equation}
    \psi(CM_{n+1})=\frac{1}{u_1 u_{1\dots 2}}\circ_2
    \psi(CM_{n})-\frac{1}{u_{1\dots 2} u_2}\circ_1 \psi(CM_{n}).
  \end{equation}
  Explicitly, $\psi(CM_{n+1})$ is given by
  \begin{equation}
    \frac{u_{2 \dots n+1}}{u_1 u_{1\dots n+1}}
    \psi(CM_{n})(u_2,\dots,u_{n+1})-\frac{u_{1\dots n}}{u_{1\dots n+1} u_{n+1}} \psi(CM_{n}).
  \end{equation}
  Then one can use the addition rule for binomial coefficients and the
  induction hypothesis.
\end{proof}

Another interesting and natural mould $PO$ has the following
components:
\begin{equation}
  PO_n=\frac{\prod_{i=2}^{n} u_{1 \dots i-1}+t\, u_i}{u_1 \prod_{i=2}^{n} (u_i u_{1\dots i})},
\end{equation}
with a parameter $t$. This mould also belongs to the image of
the Dendriform operad, as it satisfies the following equation:
\begin{equation}
 PO_{n+1}=t PO_n \succ 1/u_1 + PO_n \prec 1/u_1.  
\end{equation}
Obviously, its image by $\forg$ is given by the well-known exponential
generating series for the Stirling numbers of the first kind:
\begin{equation}
  \frac{{(1-x)}^{-t}-1}{t}.
\end{equation}

\section{Relation with quivers and tilting modules}

There is a nice relationship with the theory of tilting modules for
the equi-oriented quivers of type $\TA$ (in the classical list of
simply-laced Dynkin diagrams). Some properties of this special case
may be true in the general case of a Dynkin quiver.

\smallskip

Let $Q$ be the equi-oriented quiver of type $\TA_n$. It is known by a
theorem of Gabriel that there is a bijection between indecomposable
modules for $Q$ and positive roots for the root system of type
$\TA_n$. These positive roots are the sums $\alpha_i+\dots+\alpha_j$ for
$1\leq i\leq j\leq n$, where $\alpha_1,\dots,\alpha_n$ are the simple
roots. There is an obvious bijection $\udim$ from the set of positive
roots to the set of linear functions $u_{i \dots j}$ for $1\leq i\leq
j\leq n$, which is induced by the bijection $\alpha_i \mapsto u_i$.

A tilting module $T$ for the quiver $Q$ is a direct sum of $n$
pairwise non-isomorphic indecomposable modules such that $T$ has no
self-extension. One can therefore describe a tilting module $T$ as a
set of positive roots, satisfying some condition. Taking the inverse
of the product over the corresponding set of linear functions in the
variables $u$, one gets a rational function $\psi(T)$ for each tilting
module $T$. One can check that this set of rational functions is
exactly the image in $\mld$ of the set of planar binary trees in
$\dend$ by the operad morphism $\psi$. This gives a natural bijection
between tilting modules and planar binary trees.

By this correspondence between tilting modules and trees, the action
of the anticyclic rotation $\tau$ on the vector space $\dend(n)$ is
mapped to the action induced on the set of roots by the
Auslander-Reiten functor on the derived category of the quiver $Q$.

On the other hand, the action of the anticyclic rotation $\tau$ on the
vector space $\dend(n)$ has been related in \cite{canabull} to the
square of the Auslander-Reiten translation for the derived category of
the Tamari poset, which is a classical partial order on the set of
planar binary trees. The Tamari poset also has a very natural
interpretation in the setting of tilting modules, as a special case of
the natural partial order defined by Riedtmann and Schofield
\cite{riedscho} on the set of tilting modules of a finite-dimensional
algebra.

\smallskip

One can try and generalize this to any quiver $Q$ of finite Dynkin
type, that is any quiver whose underlying graph is a Dynkin diagram of
type $\TA$, $\TD$ or $\TE$. The theorem of Gabriel still holds, hence there
is a bijection between indecomposable modules and positive roots. One
can similarly map positive roots to linear functions in variables $u$
using the decomposition in the basis of simple roots. For an
indecomposable module $M$, the corresponding linear function is
\begin{equation}
  \udim(M)=\sum_{i=1}^{n} \dim M_i \, u_i.
\end{equation}

There is a finite set of tilting modules for $Q$. One can, just as
above, define a rational function for each tilting module $T$, as the
inverse of the product of the linear functions over the summands of
$T$. For a tilting module $T=\oplus_j M_j$, one gets
\begin{equation}
  \psi(T)=1/\prod_j \udim(M_j) 
\end{equation}

Then, one can ask the following questions:

Question 1: are the functions $\psi(T)$ for all tilting modules $T$ 
linearly independent ?

Let $V_Q$ be the vector space spanned by the $\psi(T)$ for all tilting
modules $T$.

Question 2: is $V_Q$ stable by the action induced by the action of the
Auslander-Reiten translation $\tau$ for $Q$ on the set of positive roots ?

Question 3: if so, is this action of the Auslander-Reiten translation
for $Q$ related to the Auslander-Reiten translation for the poset of
tilting modules for $Q$ defined by Riedtmann and Schofield
\cite{riedscho,happunger}.

\smallskip

Let us note a result in the same spirit, that we have learned from L.
Hille \cite{hille}: for any Dynkin quiver $Q$, one has
\begin{equation}
  \sum_T \psi(T)=1/u_1 \dots u_n,
\end{equation}
where the sum runs over the set of isomorphism classes of tilting
modules. This identity comes from a fan related to tilting modules.

\bibliographystyle{alpha}
\bibliography{moules}

\end{document}